\documentclass[twocolumn]{article}
\usepackage{authblk}
\usepackage{amsmath}
\usepackage{amssymb}
\usepackage{mathtools}
\usepackage[pdftex]{graphicx}
\usepackage[margin=20truemm]{geometry}

\author[1]{Yasuhiro Matsumoto${}^{*}$}
\affil[1]{\small Global Scientific Information and Computing Center, Tokyo Institute of Technology, Tokyo, Japan}
\affil[$*$]{\small Corresponding author; matsumoto@gsic.titech.ac.jp}
\date{}

\title{Fast wavefield evaluation method based on modified proxy-surface-accelerated interpolative decomposition for two-dimensional scattering problems}

\begin{document}

\maketitle

\abstract{This paper presents a fast wavefield evaluation method for two-dimensional wave scattering problems.
  The proposed method is based on a modified version of proxy-surface-accelerated interpolative decomposition, making it effective even if the evaluation points are near the boundary.
  The commonly known fast multipole method requires the use of direct evaluations near the boundaries of scatterers because the analytical expansion of kernel functions does not converge.
  On the one hand, the proposed method does not require the
  analytical expansion of kernel functions. The validity and effectiveness of the proposed method are demonstrated using numerical examples.
}\\

\noindent {\bf Keywords}: boundary integral equation, fast wavefield evaluation, proxy-surface method

\section{Introduction}
Electromagnetic or acoustic wave scattering problems are important in physics and engineering.
They have several interesting applications, such as invisible cloaking \cite{omni}.
The boundary integral equation method is particularly effective for the analysis of wave scattering in an infinite domain \cite{bie}.

Solving a partial differential equation using the boundary integral equations involves two steps.
In the first step, the solution on the boundaries of scatterers is obtained, and in the second step, the wavefield in the domain is evaluated using the boundary values obtained in the first step.
Although both steps have enormous computational complexity,
the fast multipole method (FMM)\cite{fmm} can quickly perform them.

The FMM, however, requires the analytical expansion of the kernel functions.
Therefore, it is not obvious how to apply the FMM to problems for which the analytical expansion of the kernel functions has not been found, such as a generalized double-layer potential that approximates the signed-distance function\cite{belyaev}.
This drawback can be overcome using algebraic techniques such as Chebyshev interpolation.
Chebyshev interpolation, however,
inherits a drawback from the analytical expansion:
it requires the use of direct evaluations near the boundaries of scatterers because the analytical expansion of kernel functions of the wave scattering problem does not converge due to a singularity at the origin of the kernel functions.
In optimization design, the average properties of the field relatively close to scatterers might need to be computed.
Therefore, a fast method that is valid even in the vicinity of boundaries is required.

To address this challenge, we propose a fast wavefield evaluation method based on a modified version of proxy-surface-accelerated interpolative decomposition (ID)\cite{fds}.
The original proxy surface method can approximate the discretized layer potential in low-rank form even for adjacent clusters.
Our method exploits this property to be fast even for a region near scatterers and to avoid the analytical expansion of kernel functions.
The proposed method thus overcomes the drawbacks of the FMM in wavefield evaluation.

The rest of this paper is organized as follows.
Section 2 formulates the wave scattering problem.
Section 3 presents the integral representations of the solution to the wave scattering problem.
Section 4 describes the proposed method in detail.
Section 5 gives the numerical results that demonstrate the validity and effectiveness of the proposed method.
Finally, Section 6 concludes this paper.

\section{Formulation of scattering problem}
Let $\Omega_2 \subset \mathbb{R}^2$ be a single (bounded) scatterer for simplicity and $\Gamma = \partial \Omega_2$ be smooth.
There is a plane incident wave $u^I$ in $\Omega_1 = \mathbb{R}^2 \setminus \overline{\Omega_2}$.
We consider the two-dimensional (2D) Helmholtz transmission problem to find solution $u$ that satisfies the following conditions:
\begin{align}
  \Delta u(x) + k_{j}^{2} u(x) = 0 \,\, {\rm in} \ \Omega_{j}, \label{eq:hel} \\
  u^{+}(x) = u^{-}(x) \, (= u(x)) \,\, {\rm on} \ \Gamma, \label{eq:bcu}\\
  \frac{1}{\varepsilon_{1}} \frac{\partial u^{+}}{\partial n}(x)
  = \frac{1}{\varepsilon_{2}} \frac{\partial u^{-}}{\partial n}(x)
  \, (= q(x)) \,\, {\rm on} \ \Gamma, \label{eq:bcq}\\
  \text{the outgoing radiation condition for } u - u^I, \label{eq:bcr}
\end{align}
where $k_j = \omega \sqrt{\varepsilon_{j} \mu_{j}}$, $\varepsilon_{j}$, and $\mu_{j}$ are the wavenumber, relative permittivity, and relative permeability in $\Omega_{i}$, respectively.
$\omega$ is the angular frequency, $n$ is the unit normal on $\Gamma$ toward $\Omega_1$, $\frac{\partial}{\partial n}$ is the normal derivative, and superscript $+$ $(-)$ represents the trace to $\Gamma$ from $\Omega_1$ $(\Omega_2)$.

The above settings imply a 2D plane electromagnetic scattering problem.
This study assumes that $\mu_{1} = \mu_{2} = 1$ and that $\varepsilon_{j}$ is constant in $\Omega_{j}$.
Note that the proposed method can also be applied to other types of boundary condition, such as the Dirichlet and Neumann
boundary conditions.

\section{Integral representations}
After the first step, in which the solution on $\Gamma$ for the boundary value problem (BVP) \eqref{eq:hel}--\eqref{eq:bcr} is found,
the solution $u$ in $\Omega_1$ or $\Omega_2$ can be evaluated using the following integral representations:
\begin{align}
  &
  \begin{multlined}
    u(x) = \int_{\Gamma} \frac{\partial G_{1}(x - y)}{\partial n} u(y) ds(y) \\
    - \int_{\Gamma} G_{1}(x - y) q(y) ds(y) + u^{I}(x) \,\,\,\, (\text{if } x \in \Omega_{1}), \label{eq:ext}
  \end{multlined} \\
  &
  \begin{multlined}
    u(x) = - \int_{\Gamma} \frac{\partial G_{2}(x - y)}{\partial n} u(y) ds(y) \\
    + \int_{\Gamma} G_{2}(x - y) q(y) ds(y) \,\,\,\, (\text{if } x \in \Omega_{2}), \label{eq:int}
  \end{multlined} 
\end{align}
where $G_{j}$ is the fundamental solution of the 2D Helmholtz equation and is given as
\begin{align}
  G_j(x-y) = \frac{\mathrm{i}}{4}H_{0}^{(1)}(k_{j} |x - y|) \label{eq:funda}
\end{align}
for $j = 1, 2$.
In \eqref{eq:funda}, $\mathrm{i}$ is the imaginary unit, $x, y \in \mathbb{R}^2$, and $H_{0}^{(1)}$ is the zeroth-order Hankel function of the first kind.
For simplicity, we only discuss \eqref{eq:ext} from here on.
However, it should be noted that the approach in Section \ref{sec:proposed} can also be applied to \eqref{eq:int} by simply switching the wavenumber and changing the signs.

We discretize the integrals of \eqref{eq:ext} by using $N$ piecewise constant elements as
\begin{align}
  \begin{multlined}
  u(x) = \sum_{j = 1}^{N} \int_{\Gamma_{j}} \frac{\partial G_{1}(x - y)}{\partial n} ds(y)u_{j} \\
  - \sum_{j = 1}^{N} \int_{\Gamma_{j}} G_{1}(x - y) ds(y) q_j + u^{I}(x), \label{eq:ext_discr}
  \end{multlined}
\end{align}
where $\Gamma_{j}$ is $j$-th boundary element ($j = 1, \ldots, N$) and $u_{j}$, $q_{j}$ are $u$, $q$ at the midpoint of $\Gamma_{j}$, respectively.
For large $N$, a fast evaluation method is required for \eqref{eq:ext_discr}.

\section{Fast wavefield evaluation method} \label{sec:proposed}
For the boundary integral equations corresponding to BVP \eqref{eq:hel}--\eqref{eq:bcr}, we use a multi-trace formulation\cite{multi}, which allows us to apply a fast algorithm\cite{fds} for finding solutions $u, q$ on $\Gamma$ without breakdown, even for transmission problems\cite{fds_trans}.
This section describes the proposed method, which is based on a modified version of proxy-surface-accelerated ID\cite{fds}.

\subsection{Interpolative decomposition}
ID is a matrix factorization method.
A matrix is decomposed into a product of basis vectors, which are a subset of the columns or rows of the matrix itself, and linear combination parameters using column-pivoted QR decomposition.
We use ID for the low-rank approximation of the discretized layer potentials in \eqref{eq:ext_discr}.
Due to space constraints, the details are omitted here; they can be found in the original paper\cite{id}.

\subsection{Modified proxy surface method}
The proxy surface method is used to evaluate a matrix to be decomposed by ID.
In the original method\cite{fds}, in the first step, the evaluated interaction between the local virtual boundary and a featured element cluster is decomposed using ID, instead of the discretized layer potential being factorized over the entire $\Gamma$,
and then the indices of the selected columns (or rows) and the linear combination parameters are obtained.
Note that the local virtual boundary is smaller than the entire $\Gamma$.

This paper proposes a modified version of the proxy surface method for the second step, in which the wavefield is evaluated using \eqref{eq:ext_discr}.
We denote the evaluation points $x$ and integration variable $y$ in \eqref{eq:ext_discr} as the $x$-side and the $y$-side, respectively.
Note that in our method, only the $y$-side is on $\Gamma$ (both sides are on $\Gamma$ in the original method).
Considering optimization design,
we want to place a lot of $x$-side points around the scatterer $\Omega_{2}$.
Since the placement of $x$-side points is arbitrary,
these points can be placed within a single unit cell and then the unit cells are aligned such that the $x$-side points are distributed around the scatterer.
Let a unit cell include $m$ evaluation points and the number of unit cells be $p_{x}$.
Each $x$-side point $x_i$ has local index $i = 1, \ldots, m$ in the unit cell.
A one-to-one correspondence to global index $i = 1, \ldots, mp_{x}$ is provided by the local index and the associated unit cell number.
Therefore, our discussion uses the local unit cell index.
Consider a virtual closed curve, called a proxy surface, surrounding a unit cell.
After evaluating the interaction between the proxy surface and the points in the unit cell, we decompose the evaluated matrix and approximate $u(x_{i})$ as
\begin{align}
  u(x_i) = \sum_{k \in S_{x}} U_{ik}u(x_{k}) \,\,\,\, (i = 1, \ldots, m), \label{eq:uprx}
\end{align}
where $S_{x}$ is the index set of selected points by ID and the number of selected points $s_{x} \ll m$, and $U \in \mathbb{C}^{m \times s_{x}}$ is the matrix of the linear combination parameters obtained by ID.
Using the unit cell allows $U$ to be reused $p_{x} - 1$ times.
  We substitute \eqref{eq:uprx} into \eqref{eq:ext_discr} and then obtain a fast method to the $x$-side points within a unit cell as
\begin{align}
  u(x_{i}) = &\sum_{k \in S_{x}} U_{ik} \Biggl( \sum_{j = 1}^{N} \int_{\Gamma_{j}} \frac{\partial G_{1}(x_{k} - y)}{\partial n} ds(y)u_{j} \nonumber \\
  &- \sum_{j = 1}^{N} \int_{\Gamma_{j}} G_{1}(x_{k} - y) ds(y) q_j \Biggr) + u^{I}(x_{i}). \label{eq:ext_uprx}
\end{align}

Similar to the procedure used for the $x$-side, we approximate the $y$-side evaluation using the proxy surface method.
Let the set of discretized $\Gamma_{j}$ $(j = 1, \ldots, N)$ be partitioned into clusters using a binary tree.
Suppose that the number of leaf clusters of the binary tree is $p_{y}$ and that each leaf has $n$ boundary elements.
The index set of the discretized $\Gamma$ included in the $t$-th leaf cluster is expressed as $C_{t}$ $(t = 1, \ldots, p_{y})$.
Consider index $t$ $(t = 1, \ldots, p_{y})$ and a virtual closed curve (proxy surface) surrounding $C_{t}$.
After evaluating the interaction between the proxy surface and $\Gamma_{j} \in {C_{t}}$, we approximate the discretized layer potentials as
\begin{align}
  &
  \begin{multlined}
    \sum_{j \in C_{t}} \int_{\Gamma_{j}} G_{1}(x - y) ds(y) q_{j} \\
    = \sum_{l \in S_{y}^{t}} \int_{\Gamma_{l}} G_{1}(x - y) ds(y) \sum_{j \in C_{t}} V^{t}_{lj} q_{j}, \label{eq:v_slp}
  \end{multlined}\\
  &
  \begin{multlined}
      \sum_{j \in C_{t}} \int_{\Gamma_{j}} \frac{\partial G_{1}(x - y)}{\partial n} ds(y)u_{j} \\
      = \sum_{l \in S_{y}^{t}} \int_{\Gamma_{l}} \frac{\partial G_{1}(x - y)}{\partial n} ds(y) \sum_{j \in C_{t}} W^{t}_{lj} u_{j}, \label{eq:w_dlp}
  \end{multlined}
\end{align}
where $S_{y}^{t}$ is an index set of selected boundary elements from $C_{t}$ obtained by ID and the number of selected elements $s_{y}^{t} \ll n$, and $V^{t}, W^{t} \in \mathbb{C}^{s_{y}^{t} \times n}$ are linear combination parameters obtained by ID for single- and double-layer potentials for $C_{t}$, respectively.
Note that unlike the original proxy surface method, in the evaluation of the interactions between clusters and virtual boundaries,
boundary elements that are a subset of adjacent clusters and are contained within the virtual boundary are not considered.
This is because the proposed method does not require the consideration of interactions between boundaries.

The resulting expression obtained using \eqref{eq:ext_uprx}, \eqref{eq:v_slp}, and \eqref{eq:w_dlp} for the proposed fast wavefield evaluation for $x_{i}$ within a specific unit cell is as follows:
\begin{align}
  &u(x_i) = \nonumber \\
  &\sum_{k \in S_{s}} U_{ik} \sum_{t = 1}^{p_{y}}  \Biggl( \sum_{l \in S_{y}^{t}} \int_{\Gamma_{l}} \frac{\partial G_{1}(x_{k} - y)}{\partial n} ds(y) \sum_{j \in C_{t}} W^{t}_{lj} u_{j} \nonumber \\
  &- \sum_{l \in S_{y}^{t}} \int_{\Gamma_{l}} G_{1}(x_{k} - y) ds(y) \sum_{j \in C_{t}} V^{t}_{lj} q_{j} \Biggr) + u^{I}(x_{i}). \label{eq:result}
\end{align}
The expression in \eqref{eq:result} can be evaluated using only matrix-vector products by considering all $x$ in the unit cell together. This evaluation is fast due to the small size of the discretized layer potentials to be evaluated.

\section{Numerical examples}
\begin{figure}[hptb]
  \centering
  \includegraphics[width=\linewidth]{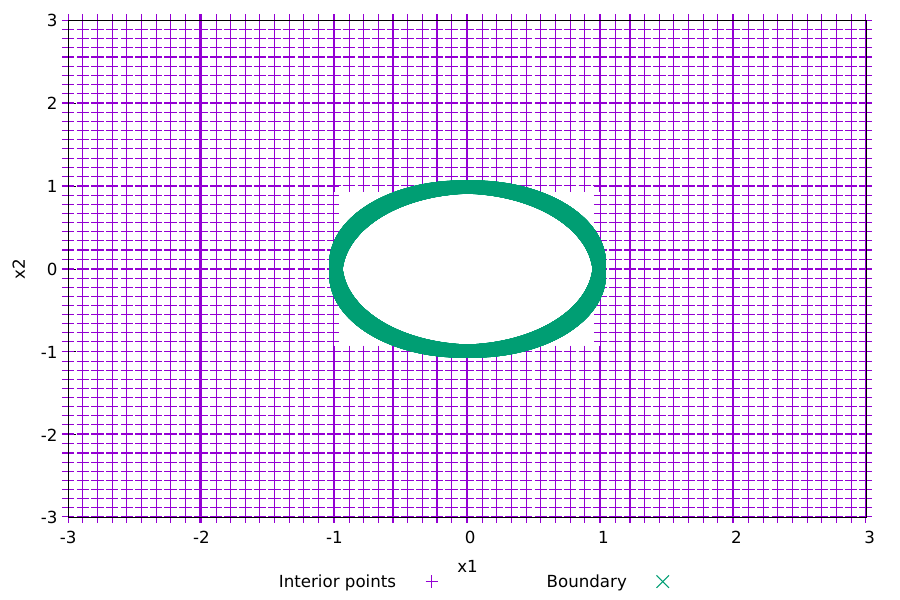}
  \caption{Evaluation points and boundary. Each evaluation point is located in a unit cell (size: $1.0 \times 1.0$). The midpoints of the discretized boundary are plotted.}
  \label{fig:x}
\end{figure}
We now demonstrate some examples of the proposed method.
In this section, $\Omega_2$ is a circle with radius 0.99, $\Gamma$ is discretized into 12800 elements,
$\omega = 10$, $\varepsilon_{1} = 1$, and $\varepsilon_{2} = 2$.
We can obtained the analytical solution corresponding to BVP \eqref{eq:hel}--\eqref{eq:bcr} since the incident wave propagating in the $x_{2}$-direction is a plane wave and $\Gamma$ is a circle.
The wavefield evaluation points are distributed at 
$x = \{ (x_{1}, x_{2}) \in \mathbb{R}^{2} \mid 1 < \|x\|_{\infty} < 3 \}$; each $x$ is located in a unit cell.
The number of $x$ points is 3200. A unit cell includes 100 points and thus 32 unit cells are used, as shown in Fig.~\ref{fig:x}.
The numerical examples were calculated on the supercomputer Pegasus at the University of Tsukuba.

We tested three versions of the proposed method and a conventional method, denoted as Conv, under the above conditions.
In the first version of the proposed method, denoted as fast-U, the proposed proxy surface method (see Section \ref{sec:proposed}) was applied to the $x$-side.
In the second version, denoted as fast-UV, the proposed proxy surface method was applied to the $x$-side, and the regular surface method was applied to the $y$-side.
In the third version, denoted as fast-UV-Vtailored, the proposed proxy surface method was applied to both the $x$-side and the $y$-side.
Note that Conv and fast-U correspond to \eqref{eq:ext_discr} and \eqref{eq:ext_uprx}, respectively.
fast-UV and fast-UV-Vtailored both correspond to \eqref{eq:result},
but their proxy methods differ.

\begin{figure}[tb]
  \centering
  \includegraphics[width=0.75\linewidth]{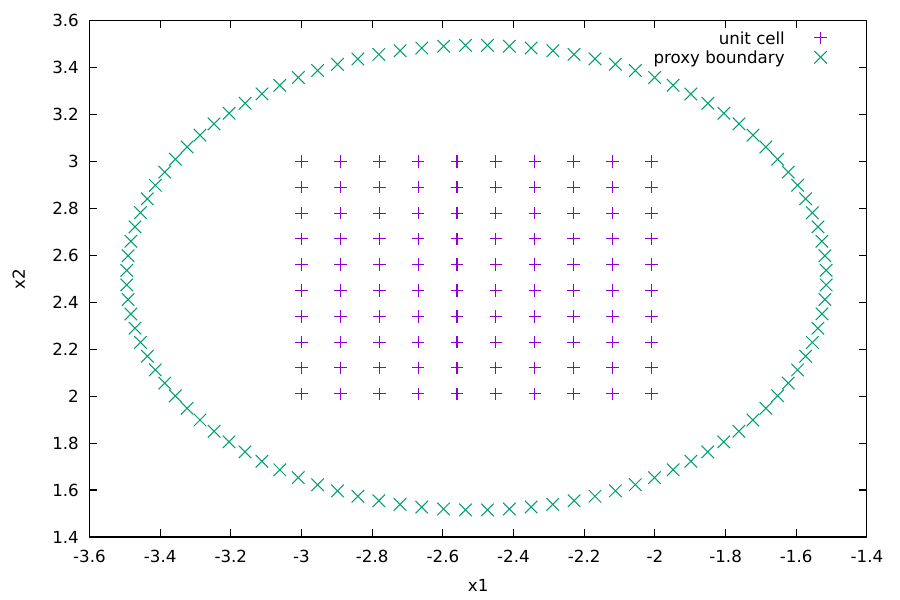}
  \includegraphics[width=0.75\linewidth]{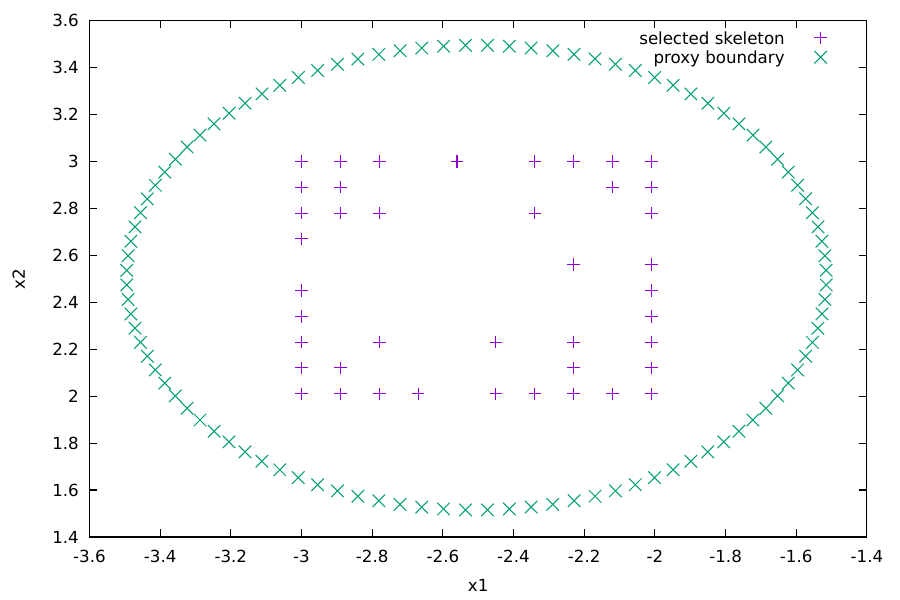}
  \caption{Results of proxy surface method for $x$-side. In this study, this result was reused in all unit cells regardless of the distance to the boundary.} Top: evaluation points in unit cell and proxy surface (local virtual boundary). Bottom: selected evaluation points (denoted as skeletons) obtained by ID and proxy surface. In this example, 42 points were selected from 100 points.
  \label{fig:proxy}
\end{figure}

Fig.~\ref{fig:proxy} shows the selected points obtained by ID with the proxy surface method applied to the $x$-side.
In all examples, we set the accuracy of ID (for the $x$-side) to $1.0 \times 10^{-6}$, which is the truncation error of column-pivoted QR decomposition for low-rank approximation.
The virtual boundary was set to a circle with a radius of 1.5 times the size of the unit cell.
The size of the unit cell was the average of the lengths in the $x_{1}$ and $x_{2}$ directions (i.e., 1).
The 100 $x$-side points in the unit cell were reduced to 42 points.
More outer points were selected because they were closer to the virtual boundary.

\begin{figure}[tb]
  \centering
  \includegraphics[width=\linewidth]{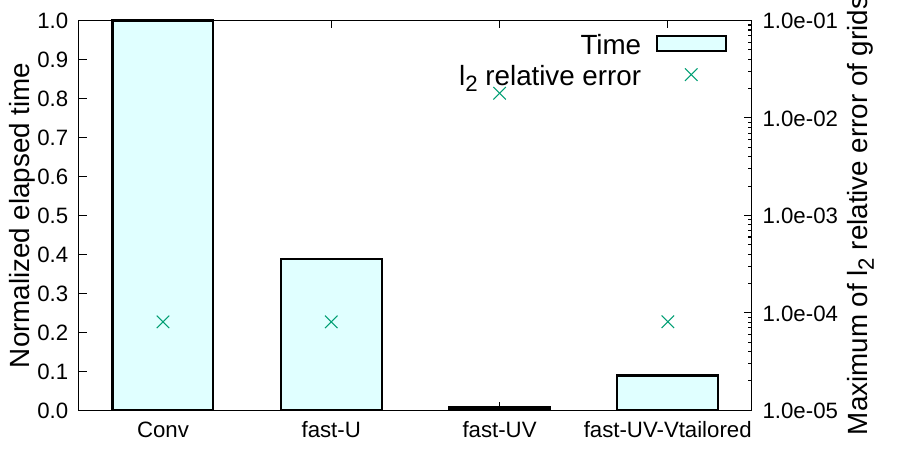}
  \caption{Comparison of normalized elapsed time (s) and $l_{2}$ relative error for various methods.
    The elapsed time of Conv was 55.1 s.
    The maximum $l_{2}$ relative errors among those evaluated for each unit cell are shown.}
  \label{fig:summary}
\end{figure}

Fig.~\ref{fig:summary} shows the improvement in computational speed for each version of the proposed method and the maximum $l_{2}$ relative error (with respect to that of the analytical solution) evaluated for each unit cell.
The $y$-side proxy method for fast-UV was a multi-level scheme\cite{fds} with an ID accuracy of $1.0 \times 10^{-12}$.
The $y$-side proxy method for fast-UV-Vtailored was a single-level scheme\cite{fds} with an ID accuracy of $1.0 \times 10^{-10}$.
Note that the difference between the single- and multi-level schemes is not essential.
The figure shows the elapsed time of each method normalized that of Conv.
The results show that fast-U, fast-UV, and fast-UV-Vtailored achieved speedups of about 2.6, 143.3, and 11.2 times, respectively, the relative errors of fast-U and fast-UV-Vtailored are almost the same as that of Conv, and fast-UV has a maximum $l_2$ relative error of 1.6\%.
If a relative error on the order of 1.0\% is acceptable, fast-UV should be selected; otherwise, fast-UV-Vtailored is a reasonable choice.
In some applications, such as optimization using a signed-distance field \cite{belyaev}, fast-UV with a relative error on the order of 1.0\% could be acceptable.

\begin{figure}[tb]
  \centering
  \includegraphics[width=0.49\linewidth]{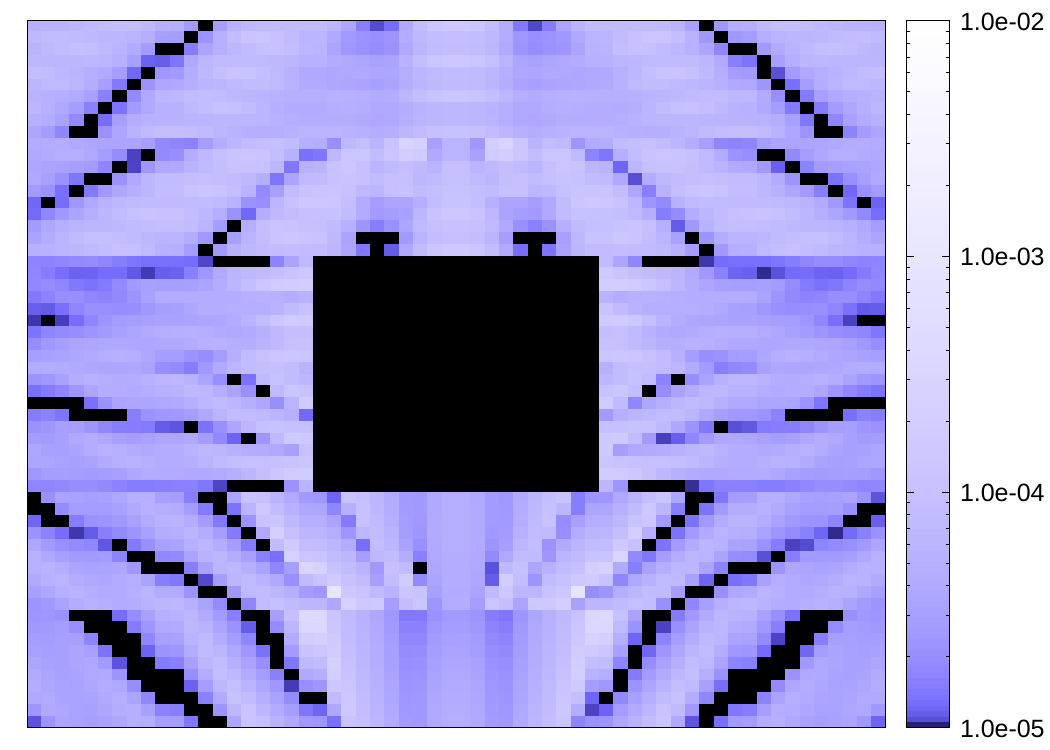}
  \includegraphics[width=0.49\linewidth]{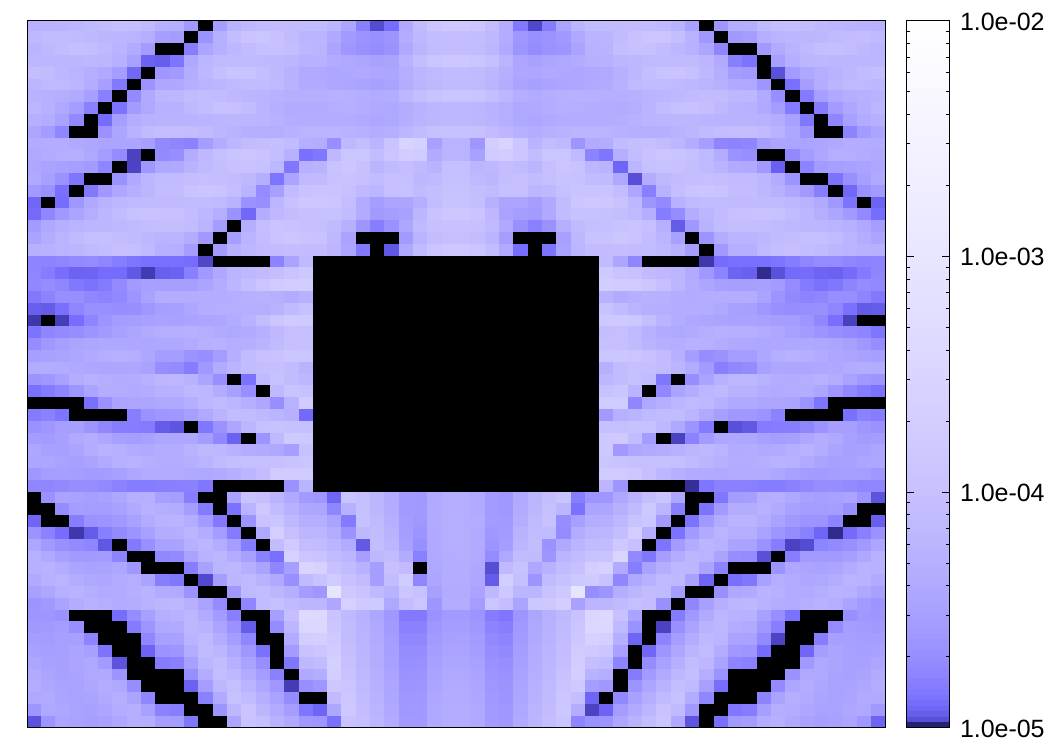}
  \includegraphics[width=0.49\linewidth]{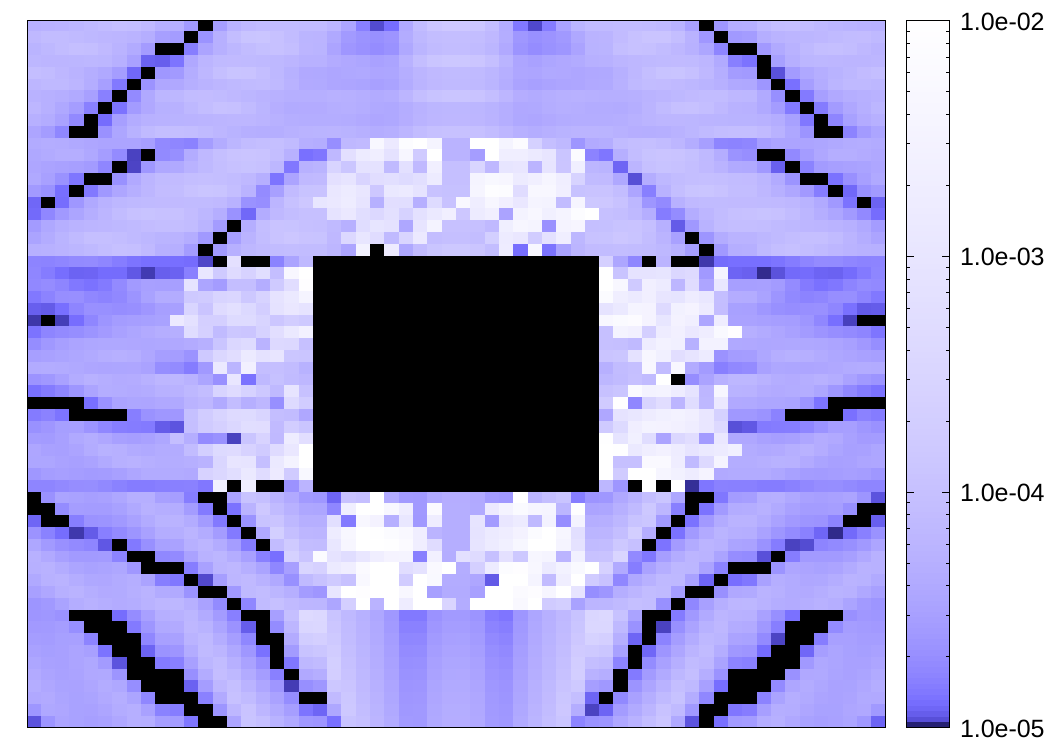}
  \includegraphics[width=0.49\linewidth]{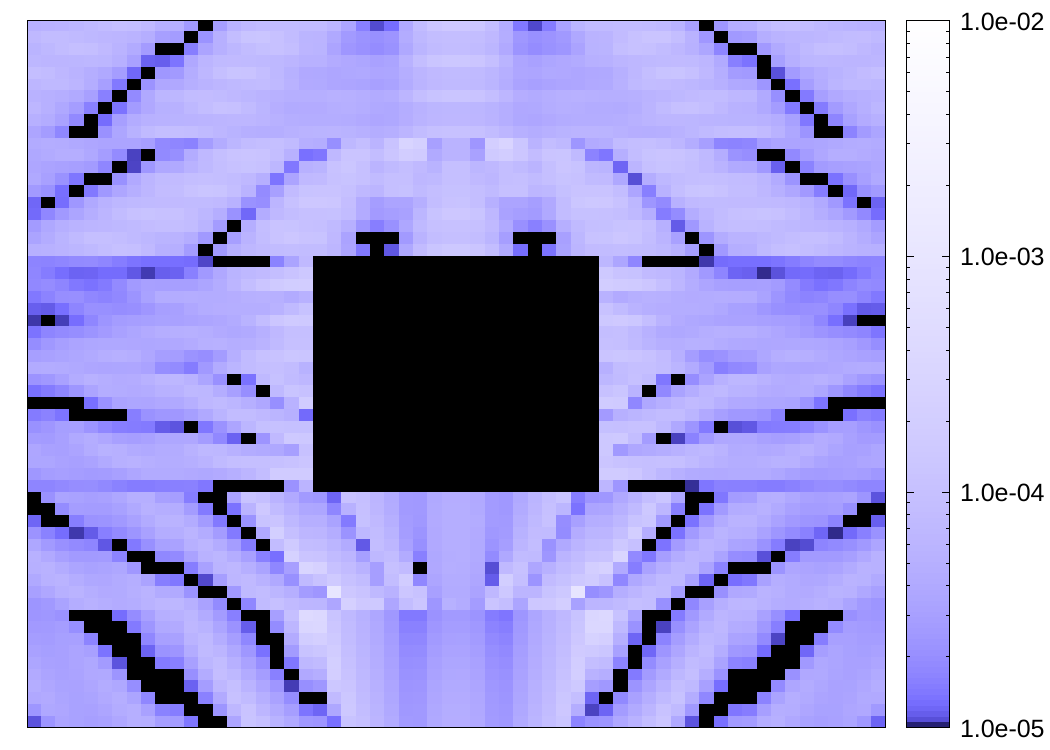}
  \caption{Relative error (color bar) of real part of total wavefield $u$ with respect to that of analytical solution at each evaluation point. The plotted domain corresponds to the region $x = \{ (x_1, x_2) \in \mathbb{R}^{2} \mid -3 < x_1 < 3, \,\, -3 < x_2 < 3 \}$. There is a scatterer in the black region. Top left: Conv. Top right: fast-U. Bottom left: fast-UV. Bottom right: fast-UV-Vtailored.}
  \label{fig:accuracy}
\end{figure}

Fig.~\ref{fig:accuracy} shows that the relative errors of the wavefield obtained by fast-U and fast-UV-Vtailored are almost the same as that obtained by Conv even near the boundary, even though the relative error for fast-UV deteriorates near the boundary.
  Note that the numerical examples in this study were calculated for a circular scatterer.
  For more complex geometries, the same proxies cannot be used in both the scatterer's nearest neighbor cells and the remaining cells to avoid loss of accuracy.

\section{Concluding remarks}
This paper proposed a fast wavefield evaluation method for 2D wave scattering.
The proposed method based on a modified version of proxy-surface-accelerated ID can be applied to a region near scatterers, which is difficult to achieve using the FMM.
The proposed method is faster than the conventional method while maintaining accuracy.

In future work, the proposed method 
will be extended to three-dimensional problems and problems for which the analytical expansion of the fundamental solution is not clear.
The optimal size of the unit cell will also be determined.


\section*{Acknowledgments}
This work was supported by the JSPS KAKENHI Grant Number 23K19972.


\begin{thebibliography}{99}

\bibitem{omni}
  K.~Matsushima, Y.~Noguchi, T.~Yamada,
  Omnidirectional acoustic cloaking against airborne sound realized by a locally resonant sonic material,
  Sci.~Rep.,
  \textbf{12.16383} (2022).

\bibitem{bie}
  J.~C.~Nédélec,
  Acoustic and electromagnetic equations: integral representations for harmonic problems,
  Springer New York,
  2001.

\bibitem{fmm}
  N.~Nishimura,
  Fast multipole accelerated boundary integral equation methods,
  Appl.~Mech.~Rev.,
  \textbf{55.4} (2002), 299--324.

\bibitem{belyaev}
  A.~Belyaev, P.~A.~Fayolle, A.~Pasko,
  Signed $L_{p}$-distance fields,
  Computer-Aided Design,
  \textbf{45.2} (2013), 523--528.

\bibitem{fds}
  P.~G.~Martinsson, V.~Rokhlin,
  A fast direct solver for boundary integral equations in two dimensions,
  J.~Comput.~Phys.,
  \textbf{205.1} (2005), 1--23.

\bibitem{multi}
  R.~Hiptmair, C.~Jerez-Hanckes,
  Multiple traces boundary integral formulation for Helmholtz transmission problems,
  Adv.~Comput.~Math.,
  \textbf{37.1} (2012), 39--91.

\bibitem{fds_trans}
  Y.~Matsumoto, N.~Nishimura,
  A fast direct solver for transmission boundary value problems for Helmholtz' equation in 2D,
  Trans.~JASCOME,
  \textbf{16} (2016), 97--102.

\bibitem{id}
  H.~Cheng, Z.~Gimbutas, P.~G.~Martinsson, V.~Rokhlin,
  On the compression of low rank matrices,
  SIAM J.~Sci.~Comput.,
  \textbf{26(4)} (2005), 1389--1404.
\end{thebibliography}
\end{document}